\documentclass[12pt]{article}
\usepackage{amsfonts}
\usepackage{graphicx}
\usepackage{amsmath}
\usepackage{amssymb}
\usepackage{psfrag}
\newcommand{\W}{\Omega}

\newcommand{\CF}{{\cal F}}

\newcommand{\CR}{\mathbb{R}}
\newcommand{\CE}{\mathbb{E}}
\newcommand{\CV}{\mathbb{V}}
\newcommand{\Ce}{{\cal E}}
\newcommand{\CP}{{\cal P}}

\title{  A law of the iterated logarithm for  sublinear expectations
\thanks{This work has been
supported by The National Basic Research Program of China (973
Program) grant No. 2007CB814901 and WCU program of the Korea Science
and Engineering Foundation (R31-20007). We thank Shige Peng and
seminar audiences at mini-workshop on $G$-expectations and
$G$-Brownian Motion at Weihai 2009, for helpful discussions and
suggestions.}}

\author{ By Zengjing  Chen,\;\;  Feng Hu
\\ Shandong University and Ajou University
}

\date{}
\begin{document}
\maketitle \noindent\textbf{Abstract}\ \ \  In this paper,
motivated  by  the notion of independent identically distributed
(IID) random variables under sub-linear expectations initiated by
Peng, we investigate a law of  the iterated logarithm for
capacities. It turns out that our theorem is a natural extension
of the Kolmogorov and  the Hartman-Wintner laws of the iterated
logarithm.




\section{ \textbf {Introduction}}

The classical laws of the iterated logarithm (LIL) as fundamental
limit theorems in probability theory play an important role in the
development of probability theory and its applications. The original
statement of the law of the iterated logarithm obtained by Khinchine
(1924) is for a class of Bernoulli random variables. Kolmogorov
(1929) and Hartman-Wintner (1941) extended Khinchine's result to
large classes of independent random variables. L\'{e}vy (1937)
extended Khinchine's result to martingales, an important class of
dependent random variables; Strassen (1964) extended
Hartman-Wintner's result to large classes of functional random
variables.  After that, the research activity of LIL has enjoyed
both a rich classical period and a modern resurgence ( see, Stout
1974 for details). To extend the LIL, a lot of fairly neat methods
have been found (see, for example, De Acosta 1983), however, the key
in the proofs of LIL is the additivity of the probabilities and the
expectations. In practice, such additivity assumption is not
feasible in many areas of applications because the uncertainty
phenomena can not be modeled using additive probabilities or
additive expectations. As an alternative to the traditional
probability/expectation, capacities or nonlinear
probabilities/expectations have been studied in many fields such as
statistics, finance and economics. In statistics, capacities have
been applied in robust statistics (Huber,1981). For example, under
the assumption of 2-alternating capacity, Huber and Strassen (1973)
have generalized the Neyman-Pearson lemma. Similarly Wasserman and
Kadane (1990) have generalized the Bayes theorem for capacities. It
is well-known that, in finance, an important question is how to
calculate the price of contingent claims. The famous Black-Shores's
formula states that, if a market is complete and self-financial,
then there exists a neutral probability measure $P$ such that the
pricing of any discounted contingent claim $\xi$ in this market is
given by $E_P[\xi].$ In this case, by Kolmogorov's strong law of
large number and LIL, one can obtain the estimates of the mean
$\mu:=E_P[\xi]$ and the variance $\sigma^2 :=E_P[|\xi-\mu|^2]$ with
probability one  by
$$
\mu=\lim_{n\to\infty} \frac 1n S_n, \quad \sigma=\limsup
\limits_{n\to\infty}(2n\log \log n)^{-1/2} |S_n-n\mu|
$$
 where $S_n$ is the sum of the first $n$ of a sample $\{X_i\}$ with
 mean $\mu$ and variance $\sigma^2.$
Statistically, an important feature of strong LLN and LIL is to
provide a frequentist  perspective for mean $\mu$ and standard
variance $\sigma.$ However, if the market is incomplete, such a
neutral probability measure is no longer unique, it is  a set $\CP$
of probability measures. In that case, one can give sub-hedge
pricing and super-hedge pricing by $\Ce[\xi]:=\inf\limits_{Q\in\CP}
E_Q[\xi]$ and $\CE[\xi]:=\sup\limits_{Q\in\CP}E_Q[\xi].$ Obviously,
both $\Ce[\cdot]$ and $\CE[\cdot]$ as functional operators of random
variables are nonlinear. Statistically, how to calculate sub-super
hedge pricing is of interest.  Motivated by sub-hedge and
super-hedge pricing and model uncertainty in finance, Peng
(2006-2009) initiated the notion of IID random variables and the
definition of $G$-normal distribution. He further obtained new
central limit theorems (CLT) under sub-linear expectations. Chen
(2009) also obtained strong laws of large numbers  in this
framework. A natural question is the following: Can the classical
LIL be generalized for capacities? In this paper, adapting the
Peng's IID notion and applying Peng's CLT under sub-linear
expectations, we investigate LIL for capacities. Our result shows
that in the nonadditive setting, the supremum limit points of $
\{(2n\log \log n)^{-1/2} |S_n|\}_{n\ge 3}$ lie, with probability
(capacity) one, between the lower and upper standard variances, the
others lie, with probability (capacity) one, between zero and the
lower standard variance. This becomes the Kolmogorov and  the
Hartman-Wintner law of the iterated logarithm if capacity is
additive, since in this case lower and upper variances coincide.

\section{ \textbf {Notations and Lemmas}}

In this section, we introduce some basic notations and lemmas. For
a given set $\CP$ of multiple prior probability measures on
$(\Omega,\CF),$ let $\mathcal {H}$ be the set of random variables
on $(\W,\CF).$

For any $\xi\in \mathcal {H},$ we define a pair of so-called
maximum-minimum expectations $(\CE,\Ce)$ by
$$\CE[\xi]:=\sup_{P\in\CP}E_P[\xi], \quad \Ce[\xi]:=\inf_{P\in\CP}
E_{P}[\xi].$$ Without confusion, here and in the sequel,
$E_P[\cdot]$ denotes the classical expectation under probability
measure $P.$ We use $\CE[\cdot]$ to denote supremum expectation
over $\CP.$

Let $\xi=I_A $ for $A\in\CF,$ immediately, a pair $(\CV,v)$ of
capacities is given by
$$\CV(A):=\sup_{P\in \CP}P(A),\;\quad  v(A):=\inf_{P\in\CP}P(A), \quad \forall A\in\CF.$$

Obviously, $\CE$ is a sub-linear expectation in the sense that\\

\noindent{\bf Definition 1} A functional $\mathbb{E}$ on $\mathcal
{H} \mapsto(-\infty,+\infty)$ is called a sub-linear expectation,
if it satisfies the following properties: for all $X,Y\in\mathcal
{H},$
\begin{description}
\item{({\bf a})} Monotonicity: $X\geq Y$ implies
$\mathbb{E}[X]\geq\mathbb{E}[Y]$.
 \item{({\bf b})} Constant preserving:
$\mathbb{E}[c]=c$, $\forall c\in\mathbb{R}$. \item{({\bf c})}
Sub-additivity: $\mathbb{E}[X+Y]\leq\mathbb{E}[X]+\mathbb{E}[Y]$.
\item{({\bf d})} Positive homogeneity: $\mathbb{E}[\lambda
X]=\lambda\mathbb{E}[X]$, $\forall\lambda\geq0$.
\end{description}

\noindent{\bf Remark} Artzner, Delbaen, Eber and Heath (1997) showed
that a sub-linear expectation indeed is a supremum expectation. That
is, if $\hat{\CE}$ is a sub-linear expectation on $\mathcal {H}$;
then there exists a set (say $\hat{\CP}$) of probability measures
such that $$\hat{\CE}[\xi]=\sup_{P\in\hat{\CP}}E_P[\xi], \quad
-\hat{\CE}[-\xi]=\inf_{P\in\hat{\CP}} E_{P}[\xi].$$ Moreover, a
sub-linear expectation $\hat{\CE}$ can generate a pair
$(\hat{\CV},\hat{v})$ of capacities denoted by
$$\hat{\CV}(A):=\hat{\CE}[I_A], \quad \hat{v}(A):=-\hat{\CE}[-I_A], \quad \forall A\in{\cal F}.$$ Therefore,
without confusion, we sometimes call the supremum expectation as the
sub-linear expectation.

It is easy to check that the pair of capacities satisfies
$$
\CV(A)+v(A^c)=1,\quad \forall A\in\CF
$$
where $A^c$ is the complement set of $A.$

For ease of exposition, in this paper, we suppose that $\CV$  and
$v$ are continuous in the sense that\\

\noindent{\bf Definition 2} A set function $V$:
$\CF\rightarrow[0,1]$ is called a continuous capacity if it
satisfies
\begin{description}
\item {({\bf 1})} $V(\phi)=0$, $V(\Omega)=1$. \item {({\bf 2})}
$V(A)\leq V(B),$ whenever $A\subset B$ and $A,B\in\CF$. \item
{({\bf 3})} $V(A_{n})\uparrow V(A)$, if $A_{n}\uparrow A$. \item
{({\bf 4})} $V(A_{n})\downarrow V(A)$, if $A_{n}\downarrow A$,
where $A_{n}, A\in \CF$.
\end{description}

Because we investigate LIL for capacity, the notion of IID random
variables under sub-linear expectations introduced by Peng is
changed slightly ( cf.\cite{peng2,peng3,peng4,peng5,peng6}).\\

\noindent{\bf Definition 3} ({\bf IID under sublinear
expectations})
\begin{description}
\item{\bf Independence:} Suppose that $Y_1,Y_2,\cdots,Y_n$ is a
sequence of random variables such that $Y_i \in\mathcal {H}.$
Random variable $Y_n$ is said to be independent of
$X:=(Y_1,\cdots,Y_{n-1})$ under $\mathbb{E}$, if for each positive
function $\varphi$ with $\varphi(X,Y_n)\in \mathcal {H}$ and
$\varphi(x,Y_n)\in\mathcal {H}$ for each $x\in \mathbb{R}^{n-1},$
we have
$$\mathbb{E}[\varphi(X,Y_n)]=\mathbb{E}[\overline{\varphi}(X)],$$
where $\overline{\varphi}(x):=\mathbb{E}[\varphi(x,Y_n)]$ and
$\overline{\varphi}(X)\in\mathcal {H}$.

\item{\bf Identical distribution:} Random variables $X$ and $Y$
are said to be identically distributed, denoted by
$X\overset{d}{=}Y$, if for each positive function $\varphi$ such
that $\varphi(X), \; \varphi(Y)\in \mathcal {H}$,
$$\mathbb{E}[\varphi(X)]=\mathbb{E}[\varphi(Y)].$$

\item{\bf IID random variables:} A sequence of random variables
$\{X_i\}_{i=1}^\infty$ is said to be IID, if
$X_i\overset{d}{=}X_1$ and $X_{i+1}$ is independent of
$Y:=(X_1,\cdots,X_i)$ for each $i\ge 1.$
\end{description}

Borel-Cantelli Lemma is still true for capacity under some
assumptions.\\

\noindent{\bf Lemma 2} Let $\{A_n,n\geq1\}$ be a sequence of
events in $\CF$ and $(\CV,v)$ be a pair of capacities generated by
sub-linear expectation $\CE$.
 \begin{description}
 \item{(1)} If $\sum\limits_{n=1}^\infty\CV(A_n)<\infty,$ then
 $\CV\left(\bigcap\limits_{n=1}^\infty\bigcup\limits_{i=n}^\infty
 A_i\right)=0.$
\item{(2)} Suppose that  $\{A_n,n\geq1\}$ are pairwise independent
with respect to $\CV$ , i.e.,
$$\CV\left(\bigcap\limits_{i=1}^\infty
A_i^c\right)=\prod_{i=1}^\infty \CV(A_i^c).$$ If
$\sum\limits_{n=1}^\infty{v}(A_n)=\infty $, then
$$v\left(\bigcap\limits_{n=1}^\infty\bigcup\limits_{i=n}^\infty A_i\right)=1.$$
\end{description}

\noindent{\bf Proof}
$$\begin{array}{lcl}
 0&\leq& \CV(\bigcap\limits_{n=1}^\infty\bigcup\limits_{i=n}^\infty
 A_i)\\
 &\leq& \CV(\bigcup\limits_{i=n}^\infty
 A_i)\\
 &\leq &\sum\limits_{i=n}^\infty\CV(A_i) \to 0,\;  \hbox{ as}\ \; n\to
 \infty.
 \end{array}
 $$
 The proof of (1) is complete.

 $$
\begin{array}{lcl}
0&\leq &1-v(\bigcap\limits_{n=1}^\infty\bigcup\limits_{i=n}^\infty A_i)\\
&=& 1-\lim\limits_{n\to \infty}v(\bigcup\limits_{i=n}^\infty A_i)\\
&=& \lim\limits_{n\to \infty}[1-v(\bigcup\limits_{i=n}^\infty A_i)]\\
& =&\lim\limits_{n\to \infty}\CV(\bigcap\limits_{i=n}^\infty A_i^c)\\
&=&\lim\limits_{n\to \infty}\prod\limits_{i=n}^\infty
\CV(A_i^c)\\
&=&\lim\limits_{n\to
\infty}\prod\limits_{i=n}^\infty(1-v(A_i))\\
&\leq&\lim\limits_{n\to \infty}\prod\limits_{i=n}^\infty{\rm
exp}(-v(A_i))\\
&=&\lim\limits_{n\to \infty}{\rm exp}(-\sum\limits_{i=n}^\infty
v(A_i))=0.
\end{array}.$$ We complete the proof of (2).\\

\noindent{\bf Definition 4} ({\bf $G$-normal distribution}, see
Definition 10 in Peng \cite{peng3}) Given a sub-linear expectation
$\CE$, a random variable $\xi\in\mathcal {H}$ with
$$\overline{\sigma}^2=\CE[\xi^2],\ \ \ \ \underline{\sigma}^2=\Ce[\xi^2]$$ is called a $G$-normal distribution,
denoted by ${\cal N}(0;[{\underline\sigma}^2,{\overline\sigma}^2])$,
if for any bounded Lipschtiz function $\phi,$ writting
$u(t,x):=\CE[\phi(x+\sqrt t \xi)],$ $(t,x)\in[0,\infty)\times \CR$,
then $u$ is the viscosity solution of PDE:
$$
\partial_t u -G(\partial^2_{xx} u)=0,\quad  u(0,x)=\phi(x),
$$
where $G(x):=\frac 12(\overline \sigma^2 x^+-\underline \sigma^2
x^-)$ and $x^+:=\max\{x,0\}$, $x^-:=(-x)^+$.\\

The following lemma can be found in Denis, Hu and Peng
\cite{Denis}.\\

\noindent{\bf Lemma 3} Suppose that $\xi$ is $G$-normal
distributed by ${\cal
N}(0;[{\underline\sigma}^2,{\overline\sigma}^2])$. Let $P$ be a
probability measure and $\phi$ be a bounded continuous function.
If $\{B_t\}_{t\geq0}$ is a $P$-Brownian motion, then
$$\CE[\phi(\xi)]=\sup_{\theta\in
\Theta}E_P\left[\phi\left(\int_0^1\theta_sdB_s\right)\right],\quad
\Ce[\phi(\xi)]=\inf_{\theta\in
\Theta}E_P\left[\phi\left(\int_0^1\theta_sdB_s\right)\right],
$$
where $$\Theta:=\left\{\{\theta_t\}_{t\ge 0}:  \theta_t \;
\hbox{is
 $\CF_t$-adaped process such that} \;\; \underline \sigma \leq
\theta_t\leq \overline \sigma\right\},$$
$$\CF_t:=\sigma\{B_s,0\leq s\leq t\}\vee{\cal N},\ \ \ \ {\cal N}\ \hbox{is the collection of P-null subsets}.$$

For the sake of completeness, the sketched proof of Lemma 3 is given in Appendix A.\\

With the notion of IID under sub-linear expectations, Peng shows the
central limit theorem under sub-linear expectations (see Theorem
5.1 in Peng \cite{peng5}).\\

\noindent{\bf Lemma 4} ({\bf Central limit theorem under
sub-linear expectations}) Let $\{X_i\}_{i=1}^\infty$ be a sequence
of IID random variables. We further assume that
$$\CE[X_1]=\Ce[X_1]=0.$$ Then the sequence $\{\overline{S}_n\}_{n=1}^\infty$
 defined by $$\overline{S}_n:=\frac{1}{\sqrt{n}}\sum\limits_{i=1}^nX_i$$ converges in law to $\xi$, i.e.,
$$\lim\limits_{n\rightarrow\infty}\CE[\varphi(\overline{S}_n)]=\CE[\varphi(\xi)],$$
for any continuous function $\varphi$ satisfying the linear growth
condition, where $\xi$ is a $G$-normal distribution.\\

\noindent{\bf Remark 1} Suppose that
$\CE[X_1^2]={\overline{\sigma}}^2$, $\overline{\sigma}>0$ and
$\varphi$ is a convex function, then, we have,
$$\CE[\varphi(\xi)]=\frac{1}{\sqrt{2\pi{\overline{\sigma}}^2}}
\int_{-\infty}^\infty\varphi(y){\rm
exp}\left(-\frac{y^2}{2{\overline{\sigma}}^2}\right){\rm d}y.$$

\section{Main results}

In this section,  we will prove the following LIL for capacities:\\

\noindent{\bf Theorem 1} Let $\{X_n\}_{n=1}^\infty$ be a sequence
of bounded IID random variables for sub-linear expectation $\CE$
with zero means and bounded variances, i.e.,
\begin{description}
\item{(A.1)} $\CE[X_1]=\Ce[X_1]=0$, \item{(A.2)}
$\CE[X_1^2]={\overline{\sigma}}^2$,
$\Ce[X_1^2]={\underline{\sigma}}^2$,  where
$0<\underline{\sigma}\leq\overline{\sigma}<\infty$.\\
\end{description}
 Denote
$S_n=\sum\limits_{i=1}^nX_i$. Then
\begin{description}

\item {(I)}$$v\left({\underline
{\sigma}}\leq\limsup\limits_{n\rightarrow\infty}\frac{S_n}{\sqrt{2n{\rm
loglog}n}} \leq{\overline{\sigma}}\right)=1.$$
\item{(II)}$$v\left({-\overline
{\sigma}}\leq\liminf\limits_{n\rightarrow\infty}\frac{S_n}{\sqrt{2n{\rm
loglog}n}} \leq-{\underline{\sigma}}\right)=1.$$ \item{(III)}
Suppose that  $C(\{x_n\})$ is the cluster set of a sequence of
 $\{x_n\}$ in $\CR,$ then
 $$
 v\left (C\left(\left\{\frac{S_n}{\sqrt{2n{\rm
loglog}n}}\right\}\right)\backslash
\left\{\limsup\limits_{n\rightarrow\infty}\frac{S_n}{\sqrt{2n{\rm
loglog}n}},\liminf\limits_{n\rightarrow\infty}\frac{S_n}{\sqrt{2n{\rm
loglog}n}}\right\}=(-\underline\sigma, \underline\sigma)\right)=1.
 $$
\end{description}

In order to prove Theorem 1, we need the following lemmas.\\

\noindent{\bf Lemma 5} Suppose $\xi$ is distributed to $G$ normal
${\cal N}(0;[{\underline\sigma}^2,{\overline\sigma}^2])$, where
$0<\underline{\sigma}\leq\overline{\sigma}<\infty$. Let $\phi$ be
a bounded continuous function. Furthermore, if $\phi$ is a
positively even function, then, for any $b\in \CR$,
$$
e^{-\frac {b^2}{2{\underline\sigma}^2}}\Ce[\phi(\xi)]\leq
\Ce[\phi(\xi-b)].
$$

\noindent{\bf Proof} Let $P$ be a probability measure,
$\{B_t\}_{t\geq0}$ be a $P$-Brownian motion. Since $\xi$ is
distributed to $G$-normal,  by Lemma 3, we have
$$
\Ce[\phi(\xi)]=\inf_{\theta\in
\Theta}E_P\left[\phi\left(\int_0^1\theta_s dB_s\right)\right].
$$

For any $\theta \in\Theta,$ write $\tilde{B_t}:=B_t-\int_0^t\frac
 b{\theta_s}ds$. By Girsanov's theorem, $\{\tilde B_t\}_{t\geq0}$ is a
 $Q$-Brownian motion  under $Q$ denoted by
$$
\frac{dQ}{dP}:=e^{-\frac 12 \int_0^1(\frac
b{\theta_s})^2ds+\int_0^1\frac b{\theta_s}dB_s}.
$$
That is
$$
\frac{dP}{dQ}=e^{-\frac 12 \int_0^1(\frac
b{\theta_s})^2ds-\int_0^1\frac b{\theta_s}d{\tilde {B}_s}}.
$$
Thus

\begin{equation}\label{IE1}
\begin{array}{lcl}
E_P\left [\phi\left(\int_0^1\theta_s dB_s-b\right)\right]
&=&E_P\left[\phi\left(\int_0^1\theta_t d(B_t-\int_0^t\frac b{\theta_s}ds)\right)\right]\\
&=&E_P\left[\phi\left(\int_0^1\theta_s d\tilde B_s\right)\right] \\
&=&E_Q\left[\phi\left(\int_0^1\theta_s d\tilde B_s\right)\cdot
e^{-\frac
12\int_0^1(\frac b{\theta_s})^2ds-\int_0^1\frac b{\theta_s}d\tilde B_s}\right]\\
&\ge &  e^{-\frac 12\left(\frac
b{\underline\sigma}\right)^2}E_Q\left[\phi\left(\int_0^1\theta_s
d\tilde B_s\right)\cdot e^{-\int_0^1\frac b{\theta_s}d\tilde
B_s}\right].
\end{array}
\end{equation}

We now prove that if $\phi$ is even, then
$$
E_Q\left[\phi\left(\int_0^1\theta_s d\tilde B_s\right)\cdot
e^{-\int_0^1\frac b{\theta_s}d\tilde B_s}\right]\ge
E_Q\left[\phi\left(\int_0^1\theta_s d\tilde B_s\right)\right].
$$
In fact, let $\overline B_t:=-\tilde B_t,$ then $\{\overline
B_t\}_{t\geq0}$ is also a $Q$-Brownian motion. Note that the
assumption that function $\phi$ is even, therefore
$$
\begin{array}{lcl}
E_Q\left[\phi\left(\int_0^1\theta_s d\tilde B_s\right)\cdot
e^{-\int_0^1\frac b{\theta_s}d\tilde B_s}\right]
&=&E_Q\left[\phi\left(\int_0^1\theta_s d\overline
B_s\right)\cdot e^{-\int_0^1\frac b{\theta_s}d\overline B_s}\right]\\
&=& E_Q\left[\phi\left(-\int_0^1\theta_s d\tilde
B_s\right)\cdot e^{\int_0^1\frac b{\theta_s}d \tilde B_s}\right]\\
 &=&
E_Q\left[\phi\left(\int_0^1\theta_s d\tilde B_s\right)\cdot
e^{\int_0^1\frac b{\theta_s}d \tilde B_s}\right].
\end{array}
$$
Since $ \frac {e^{\int_0^1b/\theta_sd\tilde B_s}+
e^{-\int_0^1b/\theta_sd\tilde B_s}}2\ge 1,$ we have

$$
\begin{array}{lcl}
E_Q\left[\phi\left(\int_0^1\theta_s d\tilde B_s\right)\cdot
e^{-\int_0^1\frac b{\theta_s}d\tilde B_s}\right] &=&\frac 12
E_Q\left[\phi\left(\int_0^1\theta_s d\tilde B_s\right)\cdot\left(
e^{\int_0^1b/\theta_sd\tilde B_s}+
e^{-\int_0^1b/\theta_sd\tilde B_s}\right)\right]\\
&\ge & E_Q\left[\phi\left(\int_0^1\theta_s d\tilde
B_s\right)\right].
\end{array}
$$

From (1), we have
$$
\begin{array}{lcl}
\Ce[\phi(\xi-b)]=\inf\limits_{\theta\in\Theta}E_P\left[\phi\left(\int_0^1\theta_s
d B_s-b\right)\right] &\ge & e^{-\frac 12(b/\underline\sigma)^2}
\inf\limits_{\theta\in\Theta} E_Q\left[\phi\left(\int_0^1\theta_s d\tilde B_s\right)\right]\\
&=& e^{-\frac
12(b/\underline\sigma)^2}\inf\limits_{\theta\in\Theta}
E_P\left[\phi\left(\int_0^1\theta_s d B_s\right)\right]\\
&=& e^{-\frac 12(b/\underline\sigma)^2} \Ce[\phi(\xi)].
\end{array}
$$

The proof of Lemma 5 is complete.\\

\noindent{\bf Lemma 6} Under the assumptions of Theorem 1, then,
for each $r>2$, there exists a positive constant $K_r$ such that
$$\CE[ \max\limits_{i\leq n}|S_{m,i}|^r]\leq
K_rn^\frac{r}{2}\ \ \ \hbox{for all }\ \ m\geq0,$$ where
$S_{m,n}=\sum\limits_{i=m+1}^{m+n}X_i$. \\

\noindent{\bf Proof.} First, we prove that there exists a positive
constant $C_r$ such that \begin{equation}\label{R3}
\sup\limits_{m\geq0}\CE|S_{m,n}|^r\leq C_rn^\frac{r}{2}.
\end{equation}

By Lemma 4 and Remark 1, it is easy to check that
$$\lim\limits_{n\rightarrow\infty}\CE[|S_{m,n}/\sqrt{n}|^r]
=\CE[|\xi|^r]=\frac{1}{\sqrt{2\pi{\overline{\sigma}}^2}}
\int_{-\infty}^\infty|y|^r{\rm
exp}\left(-\frac{y^2}{2{\overline{\sigma}}^2}\right){\rm
d}y<\infty.$$ So, we can choose
$$D_r>\frac{1}{\sqrt{2\pi{\overline{\sigma}}^2}}
\int_{-\infty}^\infty|y|^r{\rm
exp}\left(-\frac{y^2}{2{\overline{\sigma}}^2}\right){\rm d}y,$$
then there exists $n_0$ such that $\forall n\geq n_0$,
$$\CE|S_{m,n}|^r\leq D_rn^\frac{r}{2}.$$

Note that $\{X_n\}_{n=1}^\infty$ is a bounded sequence, then there
exists a constant $M>0$, such that, for each $n$, $|X_n|\leq M$. So
we can obtain (2) holds. Hence, in a manner similar to Theorem 3.7.5
of Stout \cite{stout}, we can obtain
$$\CE[ \max\limits_{i\leq n}|S_{m,i}|^r]\leq
K_rn^\frac{r}{2}\ \ \ \hbox{for all }\ \ m\geq0.$$

\noindent{\bf Lemma 7} Under the assumptions of Theorem 1, if
$$v\left(\limsup\limits_{n\rightarrow\infty}\frac{|S_n|}{\sqrt{2n{\rm
loglog}n}} \leq{\overline{\sigma}}\right)=1,$$ then, for any $b\in
\CR$ satisfying $|b|<\underline{\sigma}$,
$$v\left(\liminf\limits_{n\rightarrow\infty}\left|\frac{S_n}{\sqrt{2n{\rm
loglog}n}}-b\right|=0 \right)=1.
$$

\noindent{\bf Proof} We only need to prove that for any
$\epsilon>0$,
$$v\left(\liminf\limits_{n\rightarrow\infty}\left|\frac{S_n}{\sqrt{2n{\rm
loglog}n}}-b\right|\leq \epsilon \right)=1.$$ To do so, we only
need to prove that there exists an increasing subsequence
$\{n_k\}$ of $\{n\}$ such that
\begin{equation}\label{R3}
v\left(\bigcap_{m=1}^{\infty}\bigcup_{k=m}^{\infty}\{|S_{n_k}/\sqrt{2n_k{\rm
loglog}n_k}-b|\leq\epsilon\}\right)=1.
\end{equation} Indeed, let us choose $n_k:=k^k$ for $k\ge 1.$
For each $t>0$, write $$N_k:=[(n_{k+1}-n_k)^2t^2/2n_{k+1}{\rm
loglog}n_{k+1}],$$
$$m_k:=[2n_{k+1}{\rm loglog}n_{k+1}/t^2(n_{k+1}-n_k)],$$ $$r_k:=\sqrt{2n_{k+1}{\rm
loglog}n_{k+1}}/tm_k.$$  Since $\{X_n\}_{n=1}^\infty$ is a
sequence of IID random variables under sub-linear expectation
$\CE$, we have
\begin{equation}\label{R3}\begin{array}{lcl}
&&v\left( \left |\frac {S_{n_{k+1}}-S_{n_k}}{\sqrt{2n_{k+1}{\rm
loglog}n_{k+1}}}-b\right|\leq \epsilon\right)=
v\left(b-\epsilon\leq
\frac{S_{n_{k+1}-n_k}}{\sqrt{2n_{k+1}{\rm loglog}n_{k+1}}}\leq b+\epsilon\right)\\
&\geq& v\left(b-\epsilon/2\leq \frac{S_{N_k
m_k}}{\sqrt{2n_{k+1}{\rm loglog}n_{k+1}}}\leq
b+\epsilon/2\right)\cdot v\left(-\epsilon/2\leq
\frac{S_{n_{k+1}-n_k}-S_{N_k m_k}}{\sqrt{2n_{k+1}{\rm loglog}n_{k+1}}}\leq\epsilon/2\right)\\
&\geq& v\left(tm_k(b-\epsilon/2)\leq \frac{S_{N_k m_k}}{r_k}\leq
tm_k(b+\epsilon/2)\right)\cdot v\left(-\epsilon/2\leq
\frac{S_{n_{k+1}-n_k-N_k m_k}}{\sqrt{2n_{k+1}{\rm loglog}n_{k+1}}}\leq\epsilon/2\right)\\
&\geq& \left(v\left(bt-\epsilon t/2\leq \frac{S_{N_k}}{r_k}\leq
bt+\epsilon t/2\right)\right)^{m_k}\cdot v\left(-\epsilon/2\leq
\frac{S_{n_{k+1}-n_k-N_k m_k}}{\sqrt{2n_{k+1}{\rm
loglog}n_{k+1}}}\leq\epsilon/2\right)\\
&\geq& \left(\Ce\left[\phi(S_{N_k}/r_k
-bt)\right]\right)^{m_k}\cdot v\left(-\epsilon/2\leq
\frac{S_{n_{k+1}-n_k-N_k m_k}}{\sqrt{2n_{k+1}{\rm
loglog}n_{k+1}}}\leq\epsilon/2\right),
\end{array}
\end{equation} where $\phi(x)$ is a even function defined by
$$
\phi(x):=\left\{
\begin{array}{l}
1-e^{|x|-\epsilon t/2},\quad\quad |x|\leq \epsilon t/2;\\
0,      \quad\quad\quad\quad  \quad\;\;  |x|>\epsilon t/2.
\end{array}
\right.
$$
Note the fact that $N_k\rightarrow\infty$, as $k\rightarrow\infty$
and applying Lemmas 4 and 5,
$$
\log \Ce[\phi(S_{N_k}/r_k -bt)]\to \log \Ce[\phi(\xi -bt)]\ge
-\frac {b^2t^2}{2\underline\sigma^2} +\log \Ce[\phi(\xi)],\ \ \
\hbox{as}\ \ \ k\to\infty. $$ Thus
\begin{equation}\label{R3}\begin{array}{lcl} &&{\rm
log}\left(\Ce[\phi(S_{N_k}/r_k -bt)]\right)^{m_k}\cdot(n_{k+1}-n_k)/2n_{k+1}{\rm loglog}n_{k+1}\\
&=&(n_{k+1}-n_k)/2n_{k+1}{\rm loglog}n_{k+1}\cdot m_k{\rm
log}\Ce[\phi(S_{N_k}/r_k -bt)]\\
&\to&t^{-2}\log \Ce[\phi(\xi -bt)]\ge-\frac
{1}{2}(b/\underline{\sigma})^2+t^{-2}\log
\Ce[\phi(\xi)].\end{array}\end{equation}
 However,
\begin{equation}\label{R3}\liminf\limits_{t\rightarrow\infty}t^{-2}\log
\Ce[\phi(\xi -bt)]\geq-\frac
{1}{2}(b/\underline{\sigma})^2.\end{equation} So, from (5) and
(6), we have, for any $\delta>0$ and large enough $t$,
\begin{equation}\label{R3}\lim\limits_{k\rightarrow\infty}{\rm
log}\left(\Ce[\phi(S_{N_k}/r_k
-bt)]\right)^{m_k}\cdot(n_{k+1}-n_k)/2n_{k+1}{\rm
loglog}n_{k+1}\geq-\frac
{1}{2}(|b|/\underline{\sigma}+\delta/2)^2.\end{equation} On the
other hand, by Chebyshev's inequality,
$$\CV\left(\left|\frac{S_{n_{k+1}-n_k-N_k m_k}}{\sqrt{2n_{k+1}{\rm
loglog}n_{k+1}}}\right|>\epsilon/2\right)\leq2(n_{k+1}-n_k-{N_k
m_k){\overline{\sigma}}^2/\epsilon^2n_{k+1}{\rm
loglog}n_{k+1}}\rightarrow0,\ \ \ \hbox{as}\ \ \
k\rightarrow\infty.$$ So, as $k\rightarrow\infty$,
\begin{equation}\label{R3}\begin{array}{lcl} &&\left(n_{k+1}-n_k\right)/2n_{k+1}{\rm
loglog}n_{k+1}\cdot{\rm log}v\left(-\epsilon/2\leq
\frac{S_{n_{k+1}-n_k-N_k m_k}}{\sqrt{2n_{k+1}{\rm
loglog}n_{k+1}}}\leq\epsilon/2\right)\\
&=&(n_{k+1}-n_k)/2n_{k+1}{\rm loglog}n_{k+1}\cdot{\rm
log}\left(1-\CV\left(\left|\frac{S_{n_{k+1}-n_k-N_k
m_k}}{\sqrt{2n_{k+1}{\rm
loglog}n_{k+1}}}\right|>\epsilon/2\right)\right)\rightarrow0.\end{array}\end{equation}
Therefore, from (4), (7) and (8), we have
$$\liminf\limits_{k\rightarrow\infty}\left(n_{k+1}-n_k\right)/2n_{k+1}{\rm
loglog}n_{k+1}\cdot{\rm log}v\left( \left |\frac
{S_{n_{k+1}}-S_{n_k}}{\sqrt{2n_{k+1}{\rm
loglog}n_{k+1}}}-b\right|\leq \epsilon\right)\geq-\frac
{1}{2}(|b|/\underline{\sigma}+\delta/2)^2.$$ Now we choose
$\delta>0$ such that $|b/\underline{\sigma}|+\delta<1$. Then, for
given $\delta>0$, there exists $k_0$ such that $\forall k\geq
k_0$,
$$\begin{array}{lcl} &&v\left( \left |\frac
{S_{n_{k+1}}-S_{n_k}}{\sqrt{2n_{k+1}{\rm
loglog}n_{k+1}}}-b\right|\leq \epsilon\right)\\
&\geq&{\rm exp}\left\{-2n_{k+1}{\rm
loglog}n_{k+1}/\left(n_{k+1}-n_k\right)\cdot\left(\left(|b/\underline{\sigma}|+\delta\right)/2\right)\right\}\\
&\sim&{\rm exp}(-(|b/\underline{\sigma}|+\delta){\rm
loglog}n_{k+1}).\end{array}$$ Thus
$$\sum\limits_{k=1}^\infty v\left( \left |\frac {S_{n_{k+1}}-S_{n_k}}{\sqrt{2n_{k+1}{\rm
loglog}n_{k+1}}}-b\right|\leq \epsilon\right)=\infty.$$ Using the
second Borel-Cantelli Lemma, we have
 $$
\liminf_{k\to\infty}\left|\frac{ S_{n_k}-S_{n_{k-1}} }
{\sqrt{2n_{k}{\rm loglog}n_{k}}}-b\right|\leq \epsilon,\quad
\hbox{a.s.} \; v.
$$
But
\begin{equation}\label{ED}
\left|\frac{S_{n_k}}{\sqrt{2n_{k}{\rm loglog}n_{k}}}-b\right|\leq
\left|\frac{ S_{n_k}-S_{n_{k-1}} }{\sqrt{2n_{k}{\rm
loglog}n_{k}}}-b\right|+\frac{|S_{n_{k-1}}|}{\sqrt{2n_{k-1}{\rm
loglog}n_{k-1}}}\frac{\sqrt{2n_{k-1}{\rm
loglog}n_{k-1}}}{\sqrt{2n_{k}{\rm loglog}n_{k}}}.
\end{equation}
Note the following fact
$$
\frac{{n_{k-1}}}{n_k}\to 0,\; \hbox{as}\ \  \; k\to \infty
$$
and
$$
\limsup_{n\to\infty} |S_n|/\sqrt{2n{\rm loglog}n}\leq \overline
\sigma,\quad \hbox{a.s.} \; v.
$$
Hence, from inequality (9), for any $\epsilon >0,$
$$
\liminf_{k\to\infty}\left|\frac{S_{n_k}}{\sqrt{2n_{k}{\rm
loglog}n_{k}}}-b\right|\leq \epsilon, \hbox{ a.s.}\; v,
$$ therefore,
$$v\left(\liminf\limits_{n\rightarrow\infty}\left|\frac{S_n}{\sqrt{2n{\rm
loglog}n}}-b\right|\leq \epsilon \right)=1.
$$ Since $\epsilon$ is arbitrary, we have
$$v\left(\liminf\limits_{n\rightarrow\infty}\left|\frac{S_n}{\sqrt{2n{\rm
loglog}n}}-b\right|=0\right)=1.$$

We complete the proof of Lemma 7.\\

\noindent{\bf The Proof of Theorem 1} (I) First, we prove that
$$v\left(\limsup\limits_{n\rightarrow\infty}\frac{S_n}{\sqrt{2n{\rm
loglog}n}} \leq{\overline{\sigma}}\right)=1.$$

For each $\epsilon>0$ and $\lambda>0$, by Markov's inequality,
\begin{equation}\label{R3}
\begin{array}{lcl}
\CV\left(\frac{S_n}{\sqrt{2n{\rm
loglog}n}}>(1+\epsilon)\overline{\sigma}\right)&=&\CV\left(\frac{S_n}{\sqrt{n\overline{\sigma}^2}}
>(1+\epsilon)\sqrt{2{\rm loglog}n}\right)\\
&\leq &{\rm exp}\left(-2(1+\epsilon)^2\lambda{\rm
loglog}n\right)\CE\left[{\rm
exp}\left(\lambda\left(\frac{S_n}{\sqrt{n\overline{\sigma}^2}}\right)^2\right)\right].
\end{array}
\end{equation}
On the other hand, by Lemma 4 and Remark 1,  we have, if
$\lambda<\frac{1}{2}$,
$$\lim\limits_{n\rightarrow\infty}\CE\left[{\rm
exp}\left(\lambda\left(\frac{S_n}{\sqrt{n\overline{\sigma}^2}}\right)^2\right)\right]=\frac{1}{\sqrt{2\pi}}
\int_{-\infty}^\infty{\rm exp}(\lambda y^2){\rm
exp}\left(-\frac{y^2}{2}\right){\rm d}y<\infty.$$ Fixing $\beta>1$,
for each $\epsilon>0$, we can choose
$\lambda_\epsilon\in(0,\frac{1}{2})$ such that
$\beta=2(1+\epsilon)^2\lambda_\epsilon>1$. So, we can choose
$$C_\epsilon>\frac{1}{\sqrt{2\pi}}\int_{-\infty}^\infty{\rm exp}\left(\lambda_\epsilon
y^2\right){\rm exp}\left(-\frac{y^2}{2}\right){\rm d}y,$$ then
there exists $n_0$ such that $\forall n\geq n_0$,
\begin{equation}\label{R3} E\left[{\rm
exp}\left(\lambda_\epsilon\left(\frac{S_n}{\sqrt{n\overline{\sigma}^2}}\right)^2\right)\right]\leq
C_\epsilon.
\end{equation}
From (10) and (11), we can obtain, $\forall n\geq n_0$,
$$\CV\left(\frac{S_n}{\sqrt{2n{\rm loglog}n}}>(1+\epsilon)\overline{\sigma}\right)\leq C_\epsilon{\rm
exp}(-\beta{\rm loglog}n).$$

Choose $0<\alpha<1$ such that $\alpha\beta>1$. Let
$n_k:=[e^{k^\alpha}]$ for $k\geq1$. Then
$$\sum_{n_k\geq n_0}\CV\left(\frac{S_{n_k}}{\sqrt{2n_k{\rm loglog}{n_k}}}
>(1+\epsilon)\overline{\sigma}\right)\leq D_\epsilon\sum_{n_k\geq n_0}k^{-\alpha\beta}<\infty,$$ where $D_\epsilon$ is a positive
constant. By the first Borel-Cantelli Lemma, we can get
$$\CV\left(\bigcap\limits_{m=1}^\infty\bigcup_{k=m}^\infty\left\{\frac{S_{n_k}}{\sqrt{2n_k{\rm loglog}{n_k}}}
>(1+\epsilon)\overline{\sigma}\right\}\right)=0.$$ Also
$$v\left(\limsup\limits_{k\rightarrow\infty}\frac{S_{n_k}}{\sqrt{2n_k{\rm
loglog}n_k}} \leq(1+\epsilon){\overline{\sigma}}\right)=1.$$

Let $M_k:=\max\limits_{n_k\leq
n<n_{k+1}}\frac{|S_n-S_{n_k}|}{\sqrt{2n_k{\rm loglog}n_k}}$ for
$k\geq1$. For each $k\geq1$, $$\frac{S_n}{\sqrt{2n{\rm
loglog}n}}\leq\frac{S_{n_k}}{\sqrt{2n_k{\rm
loglog}n_k}}\frac{\sqrt{2n_k{\rm loglog}n_k}}{\sqrt{2n{\rm
loglog}n}}+M_k\frac{\sqrt{2n_k{\rm loglog}n_k}}{\sqrt{2n{\rm
loglog}n}},$$ for $n_k\leq n<n_{k+1}$. For given $\alpha$, we
choose $p>2$ such that $p(1-\alpha)\geq2$. By Lemma 6, we get
\begin{equation}\label{R3}\sum\limits_{k=1}^\infty\CE[M_k^p]\leq
K_p\sum\limits_{k=1}^\infty\frac{(n_{k+1}-n_k)^\frac{p}{2}}
{(2n_k{\rm loglog}n_k)^\frac{p}{2}}\leq
D_p^{'}\sum\limits_{k=1}^\infty k^{-\frac{p(1-\alpha)}{2}}({\rm
log}k)^{-\frac{p}{2}}<\infty,\end{equation} where $D_p^{'}$ is a
positive constant. From (12) and by Chebyshev's inequality, for
each $\epsilon>0$,
$$\sum\limits_{k=1}^\infty\CV\left(M_k>\epsilon\right)\leq\sum\limits_{k=1}^\infty\frac{\CE[M_k^p]}{\epsilon^p}<\infty.$$ Hence, by the first
Borel-Cantelli Lemma again,
$$v\left(\limsup\limits_{k\rightarrow\infty}M_k\leq\epsilon\right)=1.$$ Since $\epsilon$ is arbitrary, we have
$$v\left(\lim\limits_{k\rightarrow\infty}M_k=0\right)=1.$$ Noting
that $$\frac{\sqrt{2n_k{\rm loglog}n_k}}{\sqrt{2n_{k+1}{\rm
loglog}n_{k+1}}}\rightarrow1,\ \ \ \hbox{as}\ \ \
k\rightarrow\infty,$$ we have
$$\begin{array}{lcl}&&v\left(\limsup\limits_{n\rightarrow\infty}\frac{S_{n}}{\sqrt{2n{\rm
loglog}n}} \leq(1+\epsilon){\overline{\sigma}}\right)\\
&\geq&
v\left(\left\{\limsup\limits_{k\rightarrow\infty}\frac{S_{n_k}}{\sqrt{2n_k{\rm
loglog}n_k}}
\leq(1+\epsilon){\overline{\sigma}}\right\}\bigcap\left\{\lim\limits_{k\rightarrow\infty}M_k=0\right\}\right)\\
&=&1,\end{array}$$ which implies
\begin{equation}\label{R3} v\left(\limsup\limits_{n\rightarrow\infty}\frac{S_{n}}{\sqrt{2n{\rm
loglog}n}}\leq{\overline{\sigma}}\right)=1.\end{equation}

Similarly, considering the sequence $\{-X_n\}_{n=1}^\infty$, it
suffices to obtain
\begin{equation}\label{R3}v\left(\limsup\limits_{n\rightarrow\infty}\frac{-S_{n}}{\sqrt{2n{\rm
loglog}n}}\leq{\overline{\sigma}}\right)=1.$$ Also
$$v\left(\liminf\limits_{n\rightarrow\infty}\frac{S_{n}}{\sqrt{2n{\rm
loglog}n}}\geq-{\overline{\sigma}}\right)=1.\end{equation}

Now we prove that
$$v\left(\limsup\limits_{n\rightarrow\infty}\frac{S_n}{\sqrt{2n{\rm
loglog}n}}\ge  {\underline {\sigma}}\right)=1.$$ Indeed, from (13)
and (14), it is easy to obtain
$$v\left(\limsup\limits_{n\rightarrow\infty}\frac{|S_{n}|}{\sqrt{2n{\rm
loglog}n}}\leq{\overline{\sigma}}\right)=1.$$ For any number
$b\in(0, \underline \sigma),$ noting the fact that
$|b|<\underline\sigma$, by Lemma 7, we have
$$v\left(\liminf\limits_{n\rightarrow\infty}\left|\frac{S_n}{\sqrt{2n{\rm
loglog}n}}-b\right|=0\right)=1,$$ which implies that
\begin{equation}\label{ED}v\left(\limsup\limits_{n\rightarrow\infty}\frac{S_n}{\sqrt{2n{\rm
loglog}n}}\ge  {\underline {\sigma}}\right)=1.\end{equation}

So, from (13) and (15), we can obtain
$$v\left({\underline
{\sigma}}\leq\limsup\limits_{n\rightarrow\infty}\frac{S_n}{\sqrt{2n{\rm
loglog}n}} \leq{\overline{\sigma}}\right)=1.$$

The proof of (I) is complete.

(II) Considering the sequence $\{-X_n\}_{n=1}^\infty$, by (I), it
suffices to obtain $$v\left({\underline
{\sigma}}\leq\limsup\limits_{n\rightarrow\infty}\frac{-S_n}{\sqrt{2n{\rm
loglog}n}} \leq{\overline{\sigma}}\right)=1.$$ Thus
$$v\left({-\overline
{\sigma}}\leq\liminf\limits_{n\rightarrow\infty}\frac{S_n}{\sqrt{2n{\rm
loglog}n}} \leq-{\underline{\sigma}}\right)=1.$$

To prove (III). We only need to prove that for any number
$b\in(-\underline \sigma, \underline \sigma),$

\begin{equation}\label{ED}v\left(\liminf\limits_{n\rightarrow\infty}\left|\frac{S_n}{\sqrt{2n{\rm
loglog}n}}-b\right|=0 \right)=1.\end{equation} Noting the fact
that $|b|<\underline\sigma$, by Lemma 7, we can easily obtain
(16).

The proof of (III) is complete.

\end{document}